\documentclass[12pt]{article}
\usepackage{lscape}
\usepackage{graphicx}
\def\re#1{(\ref{#1})}

\begin{document}
\begin{titlepage}
\vskip 4 true cm
\vskip 1 true cm
\huge{\centerline{Inference about the tail of a Distribution.}}
\huge{\centerline{Improvement on the Hill Estimator }}
\vskip 0.6 true cm
\large{\centerline{ Jean Nuyts
{\small{\footnote[1]{Jean.Nuyts@umons.ac.be,
Universit\'e de Mons,
20 Place du Parc, 7000 Mons, Belgium}}}
}}
\vskip 1.3 true cm
\noindent{\bf Abstract}
\begin{quote}
{\small{
The Hill estimator is often used to infer the power behavior in tails of
experimental
distribution functions. This estimator is known to produce bad results in
certain situations which have lead to the so-called Hill horror plots. In this
brief note, we propose an improved
estimator which is simple and coherent and often 
provides an efficient remedy in
the bad situations, especially when the distribution is decreasing slowly, 
when the 
data is restricted by external cuts to lie within a finite domain, or even 
when the distribution is increasing.
}}
\end{quote}
\vfill
\end{titlepage}

\section{Introduction}

It has been advocated that self-organization, first discovered 
to dominate sand pile formation \cite{Sand}, may very well apply 
to many financial, economic, traffic control or social phenomena. 
The general outcome of these systems 
is in general an asymptotic power like behavior 
of the experimental distribution $x^{-\mu}$ in some variable $x$.
The value of the critical exponent $\mu$ is often a major prediction
of the models. 
Let us also mention models leading to 
Pareto \cite{Pareto}, Zipf \cite{Zipf1}, 
Levy-flight \cite{Levy} or Pad\'e \cite{Pade1} 
type distributions which also provide asymptotic power 
like behaviors with some definite critical exponent.

A particularly simple estimator 
of the critical exponent for a set of data measures \cite{R}
has been provided by Hill \cite{Hill}. 
It has proved to be very useful and has been widely used. See for example the KU Leuven web page \cite{BGST}. 
Following Hill,
the estimator is first applied 
to a subset of measures corresponding 
to the highest values of the variable. 
In order to increase the statistics,
this subset
is then extended to include progressively 
further experimental
values of the variable,
every set providing a value of $\mu$. 
These values are then plotted in a so-called Hill plot \cite{DHR}. 
The subset is extended 
as long as it includes values which can be considered to lie in the asymptotic domain 
where the distribution is power like.
The values obtained when they are stable 
provide the best value of the Hill estimator, 
as usually seen on the Hill plot.  
It is then obvious that the estimator usually takes 
into account all the highest experimental values 
of the variable larger than a certain value $x\geq L$.
The data points with $x<L$ do not belong to the asymptotic region
or turn out to be unsafe to use.

Various properties, aspects and generalizations of the Hill estimator 
have been discussed and studied in numerous publications. One should mention that many questions 
about the asymptotic properties,
about asymptotic normality, 
and about the volatility of the index have been addressed
(see for example \cite{DH}, \cite{CV}, \cite{dR}, \cite{HT}, \cite{GM}). 
Some authors have proposed, 
in order to improve the high volatility of the Hill plots, 
to smooth the result by averaging the Hill estimator values 
corresponding to different numbers of order statistics \cite{RS}. 
Other authors have smoothed the Hill estimator 
by convoluting the experimental random variables 
by a kernel function together with a bandwidth parameter  
\cite{CDM}. An optimal choice of the kernel and of the parameter improves the situation greatly.

However, it should be stressed that, 
in many cases, 
experimental distributions are related to discrete phenomena, 
which carry their own natural limits. 
In other cases, large values in the data may be biased, unsafe or unreachable 
for technical reasons.  
Let us give two examples.

One example of considerable importance relates to the evaluation of risk in finance, 
in particular those related to the variations of interest rates. 
Starting from their known daily variation in the past, 
one may try to evaluate the probability 
of a major overnight variation of say 3 or 4 percent in the future. 
Using a normal law, as often used by risk managers, 
the catastrophe should occur of the order of once every 10,000 years. 
Starting from the same available data but with a power law 
this may happen tree or four times per century. 
The precise evaluation of the parameter $\mu$ is obviously of paramount importance.
In data of the Federal Reserve System, 
if any day the variation is larger than a certain predefined limit the quotation is suspended. 
Hence the data is artificially cut for high $x$. 
The relevant $X_i$ must be restricted in a finite domain extending from some value $L$ 
where asymptotics begins to some value $R$ where the data is artificially cut.

Another example occurs in sociology. Suppose that some phenomena depends on the number of inhabitants 
of cities, it is clear that no data on earth today will be obtained for towns larger than say thirty millions
inhabitants or smaller than say ten inhabitants. The available data will be restricted to a finite domain. 

To conclude, in various cases, the variable is not only discrete 
but also restricted to lie within a finite domain $D_L\leq x\leq D_R$. 
Data outside the domain is either not available or unsafe to use or 
does not correspond to the asymptotic region.  
Remark that the domain may not only have a lowest (left) bound $D_L$ 
but also a highest (right) bound $D_R$. Though the lower bound $D_L$ is taken into account
by the Hill estimator,  
the possible presence of the highest bound $D_R$ 
has not been addressed by the Hill.

In this paper, we thus wish to provide a very simple 
improvement of the Hill estimator (and of the Hill plot) 
which takes into account in a perfectly symmetrical 
way of a lowest value $L$ ($D_L\leq L$) 
and of a highest value $R$ ($R\leq D_R$) 
defining a safe domain $[L,R]$ 
where the power behavior is at work, where data exists, and from which 
the critical exponent should be inferred.  
    
The problem of existing limits on the data has been subject to a much lower scrutiny
than the Hill estimator itself.
We would like to cite the work of
Beirlant and Guillou and coworkers \cite{BG}, \cite{BGDF} dealing
in particular with insurance policies where policy provisions (deduction, limits)
constrain the data. They suppose censored data
i.e. that the number of events above the constrain value is known and 
discuss the influence of censoring. 
This is not exactly the aim of our article. 

We have focused our attention on data which exist in a finite range 
only and have obtained the simplest estimator which takes into account in a symmetrical way the
lower and higher bounds of the domain.   

Finally, we would like to mention the related problems of variance 
and bias of the estimators and especially the question of their
experimental determination. For the Hill estimator, this is delicate 
and it has given rise to much research as
is testified by numerous references \cite{CGP}, \cite{GPC}, \cite{FA}, \cite{SJ}. 
Starting from the basic definition of the variance-bias parameters applied to the improved 
estimator, a detailed discussion following the paths set up for the Hill estimator
are certainly worth further work and publications, both in theory and when actual data, 
which carry their own uncertainty, are used. 
Indeed, the experimental precision of each data point 
is not always easy to determine, may in fact depend of the value of the $X_i$ itself 
and has to be taken into account precisely. On the other hand,  
in simulations, the result depends on the details of the random number generators.

\section{The Hill Estimator}

One supposes that, based on some theoretical model, 
one predicts a distribution
density $f(x)$ which behaves asymptotically as
\begin{equation}
f(x)\approx g(x)\,x^{-\mu}\ \ \ {\rm{when}} \ \ x\rightarrow \infty
\label{distribution}
\end{equation}
where $f$ and $x$ are real and the critical exponent
$\mu$ is a positive real parameter.
When $x\rightarrow \infty$, the function $g(x)$ is supposed  to be a rather
smooth function which, for large $x$ (say $x>L$), is often assumed to become
essentially a constant $\lambda$. In those cases and within
a certain margin of error, the distribution is thus approximated in the form
\begin{equation}
f(x)\approx \lambda\, x^{-\mu}\ \ \ {\rm{when}} \ \ x\rightarrow \infty
\label{limitdistribution}
\end{equation}
which is scale free.
As explained in the introduction, it has been conjectured 
by different authors that this type of power like 
distributions and in particular 
those based on self-organized critical models, 
rather than the often used Gauss like exponential forms, 
could very well
dominate certain financial, economic or social phenomena. 

Let $Y_i$ ($Y_i>0$) be a random sample obtained from experimental data for a
phenomenon which is supposed to follow a distribution law $f$ satisfying the
requirements \re{distribution} and/or \re{limitdistribution}.
The question is to draw inference on the critical
exponent $\mu$ from the random sample.

This most important question was discussed very carefully by Hill \cite{Hill}
both
from a Bayesian and a frequentist approach. He showed that, in a first approach,
both points of view lead to the same very useful answer. His recipe for
estimating $\mu$ can be outlined as follows.
\begin{itemize}
\item
Let the set $X_i$ be the set $Y_i$ reordered (reversed order statistics) in such
a way that
\begin{equation}
      X_i \geq X_j\ \ {\rm{for}}\ \ i<j
\label{orderstat}
\end{equation}
i.e. the set $X_i$ is ordered in a decreasing fashion, $X_1$ being the highest
$Y$ value.
\item
Construct the sets of numbers $H_k$, $\alpha_k$ and
$\mu_k$ for $k>2$. The set $H_k$ is defined by
\begin{equation}
H_k =\frac{1}{k}\sum_{j=1}^{k} \ln(X_j)-\ln(X_k)
\label{Hk}
\end{equation}
which is identical to the set proposed by Hill (see \cite{Hill})
\begin{equation}
\widehat H_{r} =\frac{1}{r+1}\sum_{j=1}^{r} \ln(X_j)
         -\frac{r}{r+1}\ln(X_{r+1})
\quad({\rm{with}}\quad \widehat H_r=H_{r+1})\ .
\label{Hr}
\end{equation}
The set $\alpha_k$ is defined by
\begin{equation}
\alpha_k=1/H_k
\end{equation}
and the set $\mu_k$ defined by
\begin{equation}
\mu_k=\alpha_k+1\ .
\end{equation}
\item
For quite many distributions, Hill showed that $\mu_k$ is an 
estimator of $\mu$ improving when $k$ is increased until ``it seems unwise to proceed''.
It is the point $X_k\approx L$ at the left of which the form of the actual
distribution is not anymore approximated safely enough by the equation
\re{limitdistribution}.
In other words, this occurs when the variable $x$ 
is not large enough anymore to be in the
asymptotic region of the distribution.
Both the phenomenon and
the limitation can easily be seen by
constructing specific examples.
\end{itemize}

Unfortunately, there are simple but somewhat more subtle distributions
for which the set $\mu_k$
does not converge well when $k$ is
increased. This type of situation has led to what has been called ``horror Hill
plots''.
Such horror plots can easily be constructed by considering simple examples or by
looking at page 194 of the reference \cite{EKM}. 
As argued in the introduction, there are also cases 
when the discrete experimental 
distribution carry their own domain, bounded on the left and on the right
$D_L<x<D_R$.

In this short note, we show how to remedy some of these unfortunate situations
in a very straightforward way. We construct some numerical examples to show this
explicitly. This result will be achieved by taking into account, 
not only a left boundary $L$ but also a right boundary $R$. 
It is supposed that the data cannot be trusted 
and/or is not power like at the left of $L$ or at the right of $R$.

\section{Improvement on the Hill estimator}
We first derive a simple heuristic formula and then show how to apply it to
improve the Hill estimator.

\subsection{A simple and exact formula}
Let us first suppose that $f(x)$ is exactly, for $x>0$, of a power law form
\begin{equation}
f(x)=\lambda\, x^{-\mu}
\label{powerlaw}
\end{equation}
with an arbitrary normalization constant $\lambda$. Take two arbitrary positive
numbers $0<L<R$ and define the average value $<\ln(x)>_{LR}$
of $\ln(x)$ on the interval $[L,R]$ as
\begin{equation}
<\ln(x)>_{_{LR}}
   \ =\ \frac{\int_{L}^{R} \ln(x)f(x)dx}{\int_{L}^{R} f(x)dx}\ .
\end{equation}
After some algebra, defining
\begin{equation}
\mu=\alpha+1
\label{alpha}
\end{equation}
and for later convenience
\begin{equation}
H=\frac{1}{\alpha}\ ,
\label{H}
\end{equation}
one finds, from \re{powerlaw}, the exact relation
\begin{equation}
<\ln(x)>_{_{LR}}\ =\ \frac{1}{\alpha}
     +\frac{\ln(L)L^{-\alpha}-\ln(R)R^{-\alpha}}
                       {L^{-\alpha}-R^{-\alpha}}
\label{basic}
\end{equation}
which depends on the correction function $C(\alpha,L,R)$
\begin{equation}
C(\alpha,L,R)\ =\ \frac{\ln(L)L^{-\alpha}-\ln(R)R^{-\alpha}}
                       {L^{-\alpha}-R^{-\alpha}}\ .
\label{correction}
\end{equation}
Equation \re{basic} is our basic equation which is exactly valid for an exact
power law
distribution \re{powerlaw} and approximately correct for a distribution which
satisfies \re{distribution} and/or \re{limitdistribution}.

\subsection{The basic Hill Estimator}
The basic Hill estimator is obtained when the upper limit $R$ is chosen at
infinity. Indeed,
when $\mu$ is larger than 1 and $R\rightarrow \infty$, Equation \re{basic}
reduces to
\begin{equation}
<\ln(x)>_{L\infty}= H + \ln(L)\ . 
\label{Hill}
\end{equation}

The formula \re{Hill} can easily be used to draw inference for $\alpha$ and
$\mu$ from the $k$ highest values of the random sample $X$ with the chosen order
statistics which we consider as the conditional event. Take $<\ln(x)>_{LR}$ as
the experimental average value of $\ln(x)$ from $L=X_{k}$ to the highest empirical
value $X_1$,
\begin{eqnarray}
<\ln(x)>_{L\infty}\  &\equiv&\ \frac{1}{k}\sum_{j=1}^{k} \ln(X_j)
     \label{lnaverage1}\\
\frac{1}{k}\sum_{j=1}^{k} \ln(X_j)\ &=&\  H_k+\ln(X_k) .
\label{lnaverage}
\end{eqnarray}
We find the formula for the Hill $H_k$ estimator \re{Hk} exactly.

Remark that by using a Simpson rule 
an experimental average, slightly better than \re{lnaverage1},
is obtained by
\begin{equation}
<\ln(x)>_{L\infty}\  \equiv\frac{1}{2}(X_1+X_k) 
       +\frac{1}{k}\sum_{j=2}^{k-1} \ln(X_j).
\label{lnaverage2}
\end{equation}
The difference for the Hill estimator using the slightly 
better \re{lnaverage2} 
rather than \re{lnaverage1} is usually too minute to care.

\subsection{The Improvement}

From the equation \re{basic}, we see that the Hill 
estimator can be improved by
taking into account
not only a lowest value $L=X_l$ in the reduced empirical set
$\{X_l,X_{l-1},\dots,X_r\}$ but also a highest value $R=X_r$.
Taking again a sample estimate
for $<\ln(x)>_{LR}$ in the reduced set, one
obtains the solution $\alpha_{lr}$ of
\begin{equation}
\frac{1}{l-r+1}\sum_{j=r}^{l} \ln(X_j)=
     \frac{1}{\alpha_{lr}}+
     \frac{\ln(X_{l})X_{l}^{-\alpha_{lr}}-\ln(X_{r})X_{r}^{-\alpha_{lr}}}
                       {X_{l}^{-\alpha_{lr}}-X_{r}^{-\alpha_{lr}}}
\label{fonda}
\end{equation}
as an inference for $\mu$
\begin{equation}
\mu_{lr}=1+\alpha_{lr}\ .
\label{muk}
\end{equation}

This obviously provides an easy generalization of the Hill procedure.
It makes sense when it is
known theoretically that the density function $f(x)$ behaves as
$x^{-\mu}$ times a constant $\lambda$ or times a 
slowly varying function $g(x)$ and that the random sample is secure for
$D_L\leq L\leq x\leq R \leq D_R$. It takes into account the facts that
the exact theoretical form of $f(x)$ is not known
on the left of $x=L$ and that the data points do not extend beyond
$R=X_{r}$ because of limited statistics or when the data is poorly known 
outside a finite $D_L,D_R$ domain.
In this case, the set of $X_i$ derived from the random sample $Y$ should include
the values inside the interval only and thus take the $X_i$'s
between the two limit
points ($L<X_{l}\leq X_i\leq X_{r}<R$) into account.  

The generalised Hill plot is obtained by taking first $l=r+1$, then increasing $l$ 
until ``ìt seems unwise to proceed'' and plotting the $\mu_{lr}$ as a function of $l$.
Otherwise $l$ can be increased and $r$ decreased until ``it seems unwise to proceed''.
When the data is not biased for the largest values of $X_i$, the choice $r=1$ is optimal.

It is finally worth noting the very important fact that our basic equation 
\re{fonda} is perfectly symmetrical under the exchange $X_l\leftrightarrow X_r$. 
As a result, it will be apply and provide meaningful answer 
for distribution which increase 
($\mu$ negative) rather 
than decrease with $x$, as will be seen in the examples.

Our fundamental equation \re{fonda} is transcendental and hence requires a
further treatment. The estimation of the value of $\alpha_{lr}$ inferred from
\re{fonda} can, for
example, be obtained numerically in the two following ways

\begin{description}

\item{Method (1): Direct Evaluation.}
By using standard {\it{ad\ hoc}} computer programs, the equation \re{fonda} can
be solved numerically for $\alpha_{lr}$.

\item{Method (2): Iteration.}
A simple way to achieve the same result is as follows. Define the function
$D(\alpha,L,R)$ as the derivative of the correction function 
$C(\alpha,L,R)$
with respect to $\alpha$
\begin{eqnarray}
D(\alpha,L,R)&=&\frac{\partial}{\partial\alpha} C(\alpha,L,R)
\nonumber\\
&=& \frac{R^{\alpha}L^{\alpha}(\ln(L)-\ln(R))^2}{(L^{\alpha}-R^{\alpha})^2}
\ .
\end{eqnarray}

Take the first order $\alpha^{[1]}_{lr}$ as the Hill solution
\begin{equation}
\alpha^{[1]}_{lr}=\frac{1}{\frac{\sum_{j=r}^{l} \ln(X_j)}{l-r+1}-\ln(X_l)}
\label{al1k}
\end{equation}
and define the
successive approximations $\alpha^{[p]}_{lr}, p=2,\ldots$ by iteration as
\begin{equation}
\alpha^{[p+1]}_{lr}=\alpha^{[p]}_{lr}
        \left(1
       	 	+\frac{\alpha^{[p]}_{lr}\,\frac{\sum_{j=r}^{l} \ln(X_j)}{l-r+1}
        	-\,\alpha^{[p]}_{lr}\,C(\alpha^{[p]}_{lr},X_l,X_{r})-1}
        	{{\alpha^{[p]}_{lr}}^2\,D(\alpha^{[p]}_{lr},X_l,X_{r})-1
        		}
        \right)\ .
\label{iteration}
\end{equation}
In the right hand side, the correction function $C$ and its derivative $D$ are
estimated
for the value of $\alpha$ of the preceding iteration. Then
$\alpha^{[p]}_{lr}\rightarrow \alpha_{lr}$ for $p\rightarrow \infty$. Empirically, for
$[p]$ as low as $[4]$ or $[5]$, the approximation for $\alpha_{lr}$ and hence for 
$\mu_{lr}=\alpha_{lr}{+}1$ is already very good.
\end{description}

\section{Theoretical comparison between the Hill estimator and the new estimator}

When should the new estimator obviously be preferred to the Hill estimator?
We remark that the latter (see Eq.\re{Hk}) depends critically 
on two quantities only. These are the average value of the logarithm of 
the data points and the smallest 
$X_l$ (highest $l$) included in the sample. 
The improved estimator depends on these two quantities
but also crucially on the third one: the highest attained value $X_{r}$.

Let us justify this by discussing more carefully than in the introduction
three examples of situations, pertaining to three 
completely different domains of research (sociology, economy, high energy physics), 
where both limits $X_{l}$ and $X_{r}$ should be  
taken into account.  

\begin{itemize}

\item 
Suppose that some sociological theory predicts an asymptotic power law 
for some phenomenon (crime rate, 
economic growth, power consumption)
as a function of the number of inhabitants $N_{\rm{inh}}$ in towns, 
either when this number is high or when it is small. 
Obviously, since every country has a finite number of inhabitants, 
the number of people in any city on earth is less than some number, 
actually $11.914\times 10^6$ inhabitants for Bombay (India). 
There can be no data for larger numbers. The data is cut artificially 
on this higher side essentially by    
the fact that the earth has a finite population. 
On the other end of the statistics, when $N_{\rm{inh}}$ is small, 
no aggregation of human beings is called a city 
if it contains less than say about ten inhabitants, 
depending on the rules for defining towns in different countries. 
Even more, no town will ever be defined with a fraction of one inhabitant. 
Hence the data is cut on the lower side by country depending 
administrative decisions.
  
\item
Suppose that some economic theory predicts a power law 
for the daily variation of the interest rates 
(or of the price of some share) when they are large 
(the tail of the variations). Usually, if this variation 
on some day exceeds some predefined number 
the floor has rules to suspend the quotation. 
Hence no data will ever be produced with higher variations. 
This artificial limit, which is usually not taken 
into account in the models, is included 
in the analysis of the data by the dependence 
of the new estimator in the highest value $X_{r}$.  

\item
In high energy physics, when two beams of initial particles 
are made to scatter head on at very high speed 
along a certain direction, final particles emerge 
from the scattering region at an angle with respect 
to the direction of the initial beams.
The probability distribution is often plotted as a function 
of the transverse momentum which is
related to this angle. The interaction region 
is surrounded by detectors. In the forward and in the 
backward directions there is a blind cone 
where the scattered particles cannot 
be identified among those of the intensive initial beams. 
Hence there are artificial limits both at low 
and at high transverse momenta outside which no data is 
produced because of the physics of the measuring devices. 
When the tails of the distributions of the scattered particles 
with respect to beam are studied, 
these artificial limitations have to be taken into account. 

\item

In practical situations, to obtain the improved estimator the parameters bounds $X_l$ and $X_r$,
which
limit on the left and on the right the data points,
have to be chosen in an appropriate manner. There is, a priori, no unique method for 
choosing them. For actual phenomena, related to some definite model (for example to solutions 
of a differential equation), the model itself should provide some information on the bounds. 
The onset of the asymptotic domain of the model provides an indication on $X_l$. On the other hand, 
the measure themselves may carry a natural right bound $X_r$ where the measurements become meaningless
(for example a town of more than $5\times 10^7$ inhabitants). The estimator should be evaluated for various values of
the bounds $X_l,X_r$ and educated guesses have to be applied. In general, it is convenient to simply use the highest data point as $X_r$ unless there are good reasons to reject it. In high energy physics, ``good experimentalists''
are known to apply suitable and correct cuts to their data, depending on a deep inside knowledge of their apparatus. 

\end{itemize}

\section{Tests of the new Estimator}

In order to test the new estimator, we have used 
it for some trial functions. The new estimator performs as well or even
slightly better than the original Hill estimator 
in cases where the latter is
known to be good. When applied to cases where
the Hill
estimator is known to perform more poorly, 
the improvement is often spectacular and
the new estimator converges to a much better approximate value. 
Using
the method (2) above (see Eq.\re{iteration}), 
a good value $\alpha_{lr}$ is usually obtained
after three
steps only and always after four steps 
i.e. from $\alpha_{lr}\approx\alpha^{[5]}_{lr}$.

Let us give here a few examples where we have used the following procedure. 
The examples are divided in sets. The examples (1)-(13) are related to very simple 
distributions where the best value of $\mu$ is obtained by letting $l$ become the
largest available value in the data $l=N_{\rm{rand}}$ and $r$ the lowest i.e. $r=1$. 
In the examples (14)-(17), the Hill plot is drawn for somewhat 
more complicated distributions where the deviation from a purely power behavior
shows up for small $x$. In all the cases, $r=1$ (see \re{fonda}) has been chosen.

\begin{enumerate} 
\item
We have used FORTRAN to do numerical calculations

\item
We have used the DRNGCS module of 
the International Mathematical and Statistical Libraries (IMSL) 
to generate
randomly a certain number $N_{\rm{rand}}$ of $X_i$, 
inside a predefined interval $[D_L,D_R]$ 
and following a given distribution $P$. It should be noted 
that the normalization of $P$ is inessential.
As required by DRNGCS, this normalization is chosen 
in such a way that the cumulated distribution $N(y)$ 
\begin{equation}
N(y)=\int_{D_L}^{y} P(x)dx
\label{normdist}
\end{equation}
has $N(D_R)=1$. In our examples, the distribution $P$ 
and the cumulated distribution $N$ 
have always been defined on a grid of 1000 to 10,000 points 
regularly separated between $D_L$ and $D_R$. 
This is precise enough for our purpose. 

\item
For each example (1)-(13), we give the results in a table which includes 

\begin{itemize}
\item
The trial distribution $P$, up to an arbitrary factor, 
fed to the random generator. 
\item 
The predefined limits $D_L$ and $D_R$ chosen for the 
data $D_L\leq X_i\leq D_R$ produced by the random generator 
and ordered according to Eq.\re{orderstat}.
\item
The number $N_{\rm{rand}}$ of randomly produced $X_i$.
\item
The smallest value $L=X_{N_{\rm{ran}}}$ and highest value $R=X_{1}$ 
obtained from the random generator. Obviously $D_L\leq L < R\leq D_R$.
\item
The average value $\sigma$ of the natural logarithm of the 
$X_i$ as produced by the random generator.
\begin{equation}
\sigma = \frac{1}{N_{\rm{rand}}}\sum_{i=1}^{N_{\rm{rand}}} \ln(X_i)\ .
\label{sigma}
\end{equation}
\item
The value expected for $\mu$, the exponent of the power 
decrease of the starting distribution $P$.
\item
The value produced by the Hill estimator $\mu_{\rm{Hill}}=\mu^{[1]}=\alpha^{[1]}{+}1$ 
(see Eq.\re{al1k}).
\item
The values $\mu^{[5]}$ produced by the $4^{d}$ iteration. Quite generally, 
the values $\mu^{[4]}$ produced by the $3^{d}$ iteration are already very good 
except in example (13) below where we had to proceed to 
$\mu^{[5]}$). In almost all cases (except in example (13)), it turns out 
that $\mu^{[j]}\approx \mu^{[4]}$ for $j>4$  
with at least four significant digits well within the errors of the method.

\end{itemize}

\end{enumerate}

\noindent{\bf{Example (1)}}

\noindent Let us first give an example where the old and 
the new estimators perform both well. 
We have supposed a power law \re{powerlaw} 
with $\mu=5$, $D_L=3$ and $D_R=150$. 
In this case, $D_R$ is large enough and $R$ can be essentially 
ignored in all the formulas. Following the procedure outlined 
above we obtain the line (1) of the table.
We see that in the values 
produced by the Hill estimator and by 
the new estimator are almost equal
and close to the expected value $\mu=5$. 
The new estimator is slightly better. 
The same behavior is found when $\mu$ and 
the related $D_R$ are large enough.

\vskip 0.4 cm
\noindent{\bf{Example (2)}}

\noindent The second example is almost identical to the 
first example except that the data is artificially 
reduced to the interval $D_L=3$, $D_R=4$. 
We find the line (2) of the  table.
Here we see that in this extreme case, 
where the interval allowed has been reduced drastically, 
the Hill estimator has produced a very bad 
result $\mu_{\rm{Hill}}\approx 10$ while the new estimator 
converges toward a much better value close to the expected one.

\vskip 0.4 cm
\noindent{\bf{Example (3)}}

\noindent The third example in line (3) of the table is 
identical to the second except that the number 
of $X_i$, $N_{\rm{rand}}$ is increased to 5000. 
The new estimator behaves even better 
while the Hill one is still very bad. This behavior was expected.
It should be noted that, when $D_R$ is increased from 4 to over 100, 
the two estimators tend gradually to produce compatible results. 
There is a weak dependency in the value 
chosen for $N_{\rm{rand}}$ provided that it is chosen large enough. 

\vskip 0.4 cm
\noindent{\bf{Examples (4),(5),(6)}}

\noindent When the probability distribution is chosen 
to decrease very slowly, here $P=\lambda/\sqrt{x}$, even if 
$D_R$ is chosen rather large $D_R=150$, $D_R=1500$ or $D_R=15000$, 
we see in line (4), (5) and (6) of the table that 
the new estimator is definitively better than the Hill one. 
An increase of $N_{\rm{rand}}$ and/or of $D_R$ does not alter this 
fact drastically.  

\vskip 0.4 cm
\noindent{\bf{Examples (7),(8)}}

\noindent Let us now give an example where the initial 
distribution involves a slowly varying function 
as the consequence of a Pad\'e(1,4) type distribution 
in the positive $x$ region $x>0$
\begin{equation}
P=\lambda\,\frac{1}{1+ p_2 x^2 + p^4 x^4}\ .
\label{Pade}
\end{equation}
We expect a $\mu=4$ type behavior provided 
that $x$ is chosen sufficiently large. 
The actual left limit $L$ to be used depends 
on the values of the parameters $p_i$. 
We analyze the problem in the case
of an experimental symmetric distribution 
of variation of interest rates of a Pad\'e type 
(see \cite{Pade1} and especially the results 
obtained in the first reference). 
Data which were extracted from the Federal Reserve System 
lead to $p_2=494.7$ and $p_4=4886.0$ 
for a lag of one day and a maturity of one year. 
Here, $x$ is expressed in the convenient unit of percent/year. 
The parameters $D_R$ and $D_R$ as well as $L$ and $R$ 
are given percent/year. 
It is easy to see that for the term in $p^4$ to dominate 
over the term $p^2$ in the denominator, 
a value of $x$ of about one percent in needed. 
Hence $L$ has to be chosen greater than one percent. 
We do not expect a $x^{-4}$ behavior for smaller values. 
In lines (7) and (8) of the table, 
the left boundary $D_L$ is taken to be one percent 
while the right one $D_R$, when the floor would be supposed 
to suspend the quotation, is respectively two and five percent. 

\vskip 0.4 cm
\noindent{\bf{Examples (9),(10),(11),(12) }}

\noindent Let us now turn to two examples when there is a logarithmic factor, namely 
\begin{equation}
f=\lambda\, \frac{\ln(x)}{x} 
\label{logfactor}
\end{equation}
and when there is an inverse logarithmic factor
\begin{equation}
f=\lambda\, \frac{1}{\ln(x)x}\ . 
\label{invlogfactor}
\end{equation} 
In both cases the distribution decreases rather 
slowly so that taking into account the right boundary 
has usually a major effect. 
In these cases , we have increased the value 
of $N_{\rm{rand}}$ to 5000. 
The results depend rather weakly on this choice. 
In lines (11) and (12) of the table, 
the presence of the $\ln(x)$ factor in the denominator 
leads to an effective decrease of the probability faster than $1/x$. 
Hence, the effective $\mu^{[5]}$ is expected and turns out 
to be somewhat higher than one.  
In the opposite way, in lines (9) and (10), 
when the logarithm appears in the numerator, 
the effective $\mu^{[5]}$ is expected to be somewhat 
lower than one and it is. 

\vskip 0.4 cm
\noindent{\bf{Example (13) }}

\noindent Since the treatment of the new estimator \re{fonda} is perfectly 
symmetrical under the exchange $X_{l}\leftrightarrow X_{r}$ 
of two limits $X_{r}$ and $X_{l}$, it can be applied directly without 
any change to infer the asymptotic power behavior of distributions 
which increase with $x$ and hence have a  
power behavior with negative $\mu$ up to a smooth multiplication function. 
Let us give a last example for $\mu=-3.5$, $D_L=3.0$ and $D_R=10000.0$. 
Obviously the Hill estimator produces a wrong sign for $\mu_{\rm{Hill}}$. 
The new estimator is perfect but the iterations 
had to be carried effectively to the fourth level since the starting 
$\mu^{[1]}=\mu_{\rm{Hill}}$ is so bad. 

\vskip 0.4 cm
\noindent{\bf{Example (14). Figure 1 }}

\noindent In figure 1, we show the original Hill plot and the improved plot related to the Pad\'e example (8).
The parameters in Eq. \ref{Pade} are again $p_2=494.7$ and $p_4=4886.0$. Here we have chosen $D_L=1$ percent and 
$D_R=3$ percent and generated 2000 data points. The continuous line is the original Hill plot, 
the dotted line is the improved plot and the dashed line is the expected asymptotic value $\mu=4$.

\vskip 0.4 cm
\noindent{\bf{Example (15). Figure 2  }}

\noindent We have taken a distribution 
\begin{equation}
f(x)=3.0\ x^{-4.0}+1.0\ x^{-2.5}\ .
\end{equation}
The random data with $N_{\rm{rand}}=10000$ has been collected 
in the interval $D_L=10$ and $D_R=30$.
The original Hill estimate and the improved 
estimate are plotted in
terms of $l$ ranging up to $10000$ in figure 2. 
The continuous line is the original Hill plot, 
the dotted line is the improved plot and the dashed line is the expected asymptotic value $\mu=2.5$.

\vskip 0.4 cm
\noindent{\bf{Example (16). Figure 3  }}

\noindent In figure 3, the distribution is taken as $\ln(x)/x)$ related to example (10).
The random data with $N_{\rm{rand}}=10000$ has been collected 
in the interval $D_L=100$ and $D_R=400$. The continuous line is the original Hill plot, 
the dotted line is the improved plot and the dashed line 
is the expected asymptotic value $\mu=1$ (regardless of the $\ln(x)$). 
In fact one sees that when $k$ reaches about 4000, the value of $\mu$ is close to 1. 
For larger values of $k$, $\mu$  decreases below $\mu=1$. This is due to the fact that the 
slowly varying function $\ln(x)$ in the numerator leads 
to a slightly decreased value of $\mu$ whose effects increases with $k$. Indeed,
since $\ln(x)$ is an increasing function of $x$, the distribution is expected to decrease somewhat more slowly
than $1/x$.    

\vskip 0.4 cm
\noindent{\bf{Example (17). Figure 4  }}

\noindent In figure 4, the distribution is taken as $1/\sqrt{x}$ related to examples 
(4), (5), (6).
The random data with $N_{\rm{rand}}=10000$ has been collected 
in the interval $D_L=3$ and $D_R=1500$. The continuous line is the original Hill plot, 
the dotted line is the improved plot and the dashed line 
is the expected asymptotic value $\mu=0.5$. 

\begin{figure}[th]\label{fig1}
\caption{A Pad\'e $[0,4]$ behavior (see Example 14) }
\begin{center}
\includegraphics*[width=8.0cm]{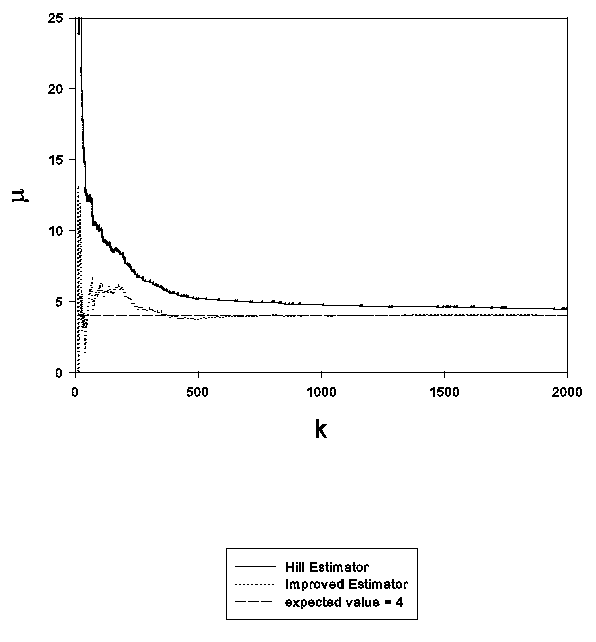}
\end{center}
\end{figure}

\begin{figure}[ht]\label{fig2}

\begin{center}\caption{A sum of two power like behaviors (see Example 15)}
\includegraphics*[width=12.0cm]{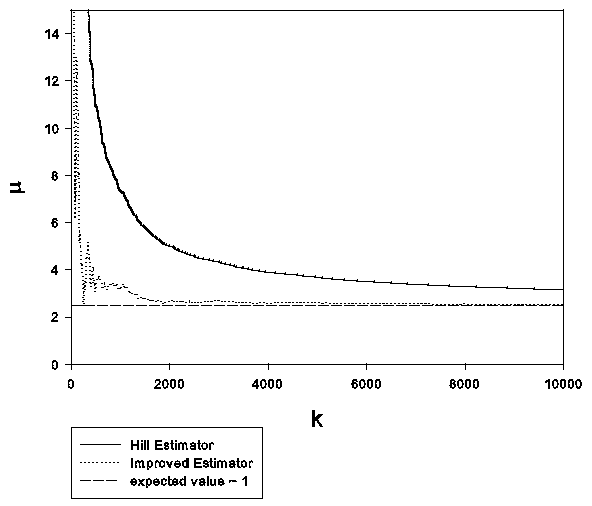}
\end{center}
\end{figure}

\begin{figure}[ht]\label{fig3}
\caption{A $\ln(x)/x$ behavior (see Example 16)}
\begin{center}
\includegraphics*[width=8.0cm]{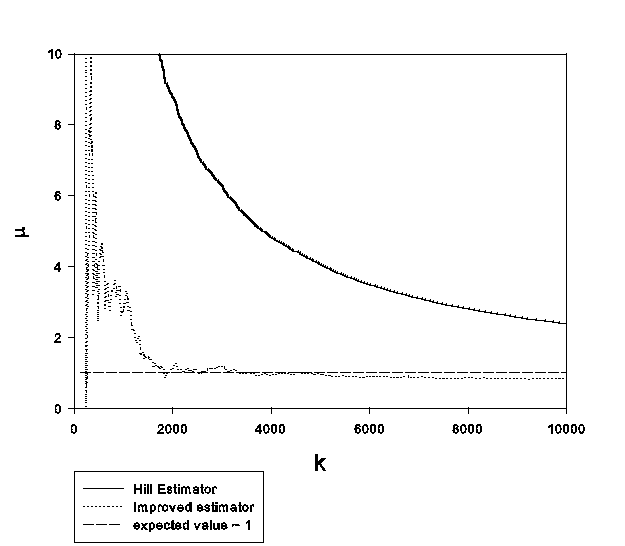}
\end{center}
\end{figure}

\begin{figure}[ht]\label{fig4}
\caption{A $1/\sqrt{x}$ behavior (see Example 17)}
\begin{center}
\includegraphics*[width=8.0cm]{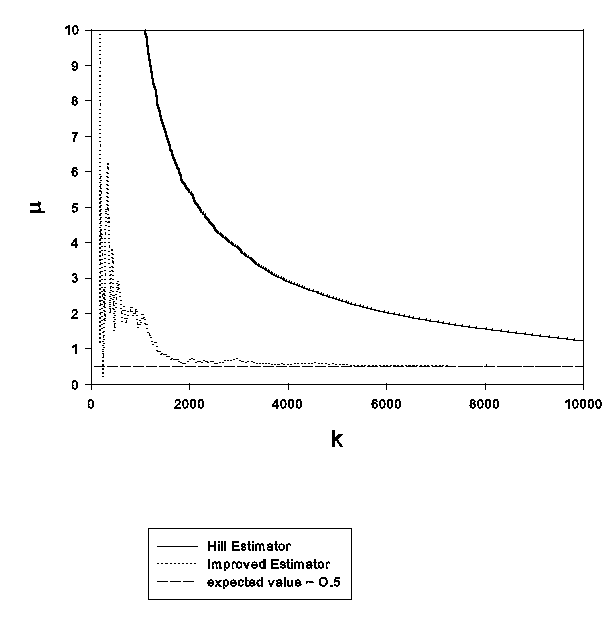}
\end{center}
\end{figure}

\begin{landscape}
\hspace*{1 cm}
\vspace*{1 cm} 
\begin{tabular}{|c|c|c|c|c|c|c|c|c|c|c|}
\hline
&$P$ & $N_{\rm{rand}}$&$D_L$&$D_R$& $L$ & $R$ &$\sigma$ & input $\mu $           
     &$\mu_{\rm{Hill}}$&$\mu^{[5]}$
      \\ \hline
(1)&$x^{-5}$&1000&3.0&150.0&3.0004&14.37&1.339&5&
 5.157&5.115
      \\ \hline
(2)&$x^{-5}$&1000&3.0&4.0&3.0015&3.962&1.208&5&
 10.175&5.784 
      \\ \hline
(3)&$x^{-5}$&5000&3.0&4.0&3.0002&3.999&1.214&5&
 9.634&5.142 
      \\ \hline
(4)&$\frac{1}{\sqrt{x}}$&1000&3.0&150.0&3.017&149.53&3.642&0.5&
 1.394&0.511 
      \\ \hline
(5)&$\frac{1}{\sqrt{x}}$&1000&3.0&1500.0&3.058&1494.7&5.585&0.5&
 1.223&0.508 
      \\ \hline
(6)&$\frac{1}{\sqrt{x}}$&1000&3.0&15000.0&3.115&14945&7.682&0.5&
 1.153&0.518 
      \\ \hline
(7)&$\frac{1}{1+p_2x^2+p_4x^4}$&1000&1.0&2.0&1.0001&1.992&0.234&$\approx$ 4&
 5.269&3.977 
      \\ \hline
(8)&$\frac{1}{1+p_2x^2+p_4x^4}$&1000&1.0&5.0&1.0001&4.686&0.320&$\approx$ 4&
 4.126&3.980 
      \\ \hline
(9)&$\frac{\ln(x)}{x}$&1000&100.0&400.0&100.07&399.11&5.322&$\approx$ 1&
 2.397&0.849 
      \\ \hline
(10)&$\frac{\ln(x)}{x}$&1000&2000.0&5000.0&2001.7&9973.2&8.423&$\approx$ 1&
 2.217&0.914 
      \\ \hline
(11)&$\frac{1}{\ln(x)x}$&5000&8000.0&10000.0&8000.8&9999.0&9.098&$\approx$ 1&
 10.054&1.247 
      \\ \hline
(12)&$\frac{1}{\ln(x)x}$&5000&3000.0&6000.0&3000.9&5998.1&8.346&$\approx$ 1&
 3.944&1.164 
      \\ \hline
(13)&$x^{3.5}$&1000&3.0&10000.0&1828.3&9995.9&8.989&-3.5&
 1.677&-3.503 
      \\ \hline     
\end{tabular}
\vskip 0.5 cm
\hskip 2 cm {\bf{Table}} (see text for the notation)
\end{landscape}

\newpage

\section{Conclusion}

We have shown that, by a very simple alteration, the Hill estimator, which
provides a useful inference for the power behavior of the tail of a
distribution, can be
easily improved to cover cases where it performs badly.
This includes cases
when the tail is not of the form of a power law but
is multiplied by a slowly varying functions $g(x)$ including logarithms for
example. It also applies when the power $\mu$ is of 
the order or smaller than one.
It even works for inferring the tails of increasing 
distributions i.e. when $\mu$ is negative. 

Needless to say, the approach outlined in this paper is one of efficiency. From
the frequentist point of view, it is perfectly justified. However, in some way,
the Bayesian point of view is also, but granted not completely, met as the
theoretical and empirical knowledge of extreme
tails of the distribution is somehow better taken into account.
It should be stressed that the new estimator applies particularly 
when the data is confined by artificial 
external conditions to lie within a finite domain bounded not 
only by a lowest value of the variable but also by a highest value.

\vfill\eject
\vskip .7 true cm
\noindent{\large\bf{Acknowledgment}}
\vskip 0.2 true cm
\noindent

This work was supported in part by
the Belgian F.N.R.S. (Fonds National de la Recherche Scientifique).
The author would like to thank Professor Isabelle Platten for carrying detailed
numerical tests in an early phase of this research and Professor Fernand Grard 
for a critical reading of the manuscript.

\goodbreak
\end{document}